\magnification=1200
\noindent


\bigskip

\centerline{\bf Number of points of Prym varieties over finite fields}

\medskip

\centerline{
Marc Perret}

\bigskip

{\it Abstract. We prove some upper and lower bounds for the number of rational points of Prym varieties defined over  finite fields.  They are better than the usual Weil bounds valid for any abelian varieties defined over such fields. 

}

\bigskip


\bigskip

\noindent {\bf I. Introduction.} Let $\pi : Y \rightarrow X$ be a covering of smooth algebraic irreducible projective curves defined over a field $k$ of zero  or odd characteristic. Then the Jacobian $J_X$ of $X$ is isogenous to a sub-abelian variety of the Jacobian $J_Y$ of $Y$. If we suppose moreover that $\pi$ have degree 2, then the non-trivial involution $\sigma$ of this covering induces an involution $\sigma^*$ on $J_Y$. The well known proposition part of the following is easy to prove :

\medskip

 \noindent {\bf Definition and proposition 1.} {\it The Prym variety $Pr = Pr_\pi$ associated to the unramified double cover $\pi : Y \rightarrow X$ of a curve $X$ of genus $g\geq 2$ is defined as $Pr := Im(\sigma^*-id)$. It is an abelian subvariety of $J_Y$ of dimension $g-1$, isogenous to a direct factor of $J_X$ in $J_Y$.}
 
 \bigskip
 
For more details, see [1] or [4]. The computation of the genus follows from Riemann-Hurwitz Theorem. It is known that Prym varieties are in general {\it not} jacobian varieties. For instance, it has been proved by Beauville in [1] that any abelian variety of dimension less than $5$ is a degeneration of Prym varieties-at least on algebraically closed fields.

Suppose from now on that $k$ is the finite field ${\bf F}_q$ with $q$ elements. Being an abelian variety of dimension $g-1$, we can apply to $Pr$ the following theorem (see the historical source [6] for instance) :

\medskip

\noindent {\bf Theorem (Weil, 1948).} {\it Let $A$ be an abelian variety of dimension $d$ defined over ${\bf F}_q$. Then there exists $\theta_1, \dots, \theta_d \in {\bf R}/2\pi{\bf Z}$, such that for any $n \geq 1$, the number of rational points of $A$ over ${\bf F}_{q^n}$ is given by
$$\sharp A({\bf F}_{q^n}) = \prod_{i=1}^d (q^n+1-2\sqrt {q^n} \cos n\theta_i). \leqno (\imath)$$
\noindent In particular,
$$(q+1-2\sqrt q)^d \leq \sharp A({\bf F}_q)\leq (q+1+2\sqrt q)^d. \leqno (\imath \imath)$$
$(\imath \imath \imath)$ If in addition $A$ is the jacobian of a curve $C$ of genus $g$, then $d=g$
and the $\theta_i$'s are also related to the number of rational points of $C$ over ${\bf F}_{q^n}$ by
$$\sharp C({\bf F}_{q^n}) = q^n+1-2\sqrt {q^n} \left( \sum_{i=1}^g \cos n\theta_i\right).$$
}

\bigskip

The $(\imath \imath)$ part of Weil theorem for the prymian variety $Pr_\pi$ of a double unramified cover $\pi$ of a curve $X$ of genus $g$ reads :

$$(q+1-2\sqrt q)^{g-1} \leq \sharp Pr_\pi({\bf F}_q) \leq (q+1+2\sqrt q)^{g-1}. \leqno (1)$$

The upper and lower bounds in $(1)$ are ``best possible" in the sense that both can be reached.
Indeed, it is known that any elliptic curve is a Prymian variety.
Now, suppose that $E$ is chosen so that it reaches the upper bound (resp. the lower bound) of Weil's inequality $(\imath \imath)$.
Such an elliptic curve does exist if $q$ is a square, see [5]. Then $E$ reaches the upper (resp. the lower) bound of $(1)$.

The existence of such an elliptic curve $E$ proves of course also that the $(\imath \imath)$ part of Weil theorem for the jacobian variety $J_X$ of a curve $X$ :
$$(q+1-2\sqrt q)^{g} \leq \sharp J_X({\bf F}_q) \leq (q+1+2\sqrt q)^{g}. \leqno (2)$$
is also ``best possible", at least for $g=1$. Nevertheless, several sharper {\it lower} bounds and an upper bound where proved in [3]  for this jacobian. For instance :

\medskip

\noindent {\bf Theorem (G. Lachaud, M. Martin-Deschamps, 1991).} {\it Let $J_X$ be the jacobian variety of a genus $g$ curve $X$ defined over ${\bf F}_q$, and $\sharp X({\bf F}_q)$ be the number of rational points of $X$. Then
$$(\sqrt q-1)^2 {q^{g-1}-1 \over g} {\sharp X({\bf F}_q) + q-1 \over q-1} \leq \sharp J_X({\bf F}_q). \leqno (3)$$
If $X$ admits a map of degree $d$ onto the projective line, then one have also
$$\sharp J_X({\bf F}_q) \leq {e \over q}  (2g\sqrt e)^{d-1} q^g. \leqno (3~bis)$$

}

\bigskip

The aim of this paper is to prove some Lachaud-Martin Deschamps type bounds for prymian varieties (see theorem 2). The method used is different from their one. It also gives some bounds for jacobians (theorem 5), but they are not always as good as $(3)$ and $(3~bis)$, see remark 3 below.

\bigskip

\noindent {\bf II. Bounds for prymian varieties.} If $C$ is an algebraic curve defined over a finite field $k$ with $q$ elements, we denote by $\sharp C({\bf F}_q)$ the number of ${\bf F}_q$-rational points of $C$. 

\medskip

The main result of this paper is the following :

\medskip

 \noindent {\bf Theorem 2.} {\it Let $X$ be an absolutely irreducible projective smooth algebraic curve defined over the finite field $k$ of odd characteristic with $q$ elements. Let $g$ be the genus of $X$, and let $\pi : Y \rightarrow X$ be an unramified covering of degree $2$. Then :

$$
\left({\sqrt q+1 \over \sqrt q-1}\right)^{{\sharp Y({\bf F}_q)-\sharp X({\bf F}_q) \over 2 \sqrt q}-2\delta} (q-1)^{g-1}
\leq \sharp Pr({\bf F}_q)\leqno (\imath)$$
with 
$\delta = 1$ if ${\sharp Y({\bf F}_q)-\sharp X({\bf F}_q)\over 2 \sqrt q} \notin {\bf Z}$, and $\delta = 0$ otherwise.

$$\sharp Pr({\bf F}_q) \leq \left(q+1 + {\sharp Y({\bf F}_q)-\sharp X({\bf F}_q)\over g-1} \right)^{g-1}.\leqno (\imath \imath)$$
$(\imath \imath \imath)$ If $X$ admits a degree $d$ map onto ${\bf P}^1_{{\bf F}_q}$, then
$$\left({\sqrt q-1 \over \sqrt q+1}\right)^{d{q+1 \over 2 \sqrt q} +2} (q-1)^{g-1}\leq \sharp Pr({\bf F}_q) \leq e^d (q+1)^{g-1}.$$

}

\bigskip

\noindent {\it Proof of theorem 2.} We use the following lemmas 3 and 4, whose proofs are postponed to the end of this section.

\medskip

 \noindent {\bf Lemma 3.} {\it
Let $a > 0$, $\gamma \in {\bf N}$, $b \in {\bf R}$ and $P$ be the polytope
$P = \{(x_1, \dots, x_\gamma) \in [-1, 1]^\gamma \vert \sum_{k=1}^\gamma x_k = b\}$.
Suppose that $\vert b \vert \leq \gamma$. Then 
$$\eqalign 
{\inf_{(x_1, \dots, x_\gamma) \in P} \prod_{k=1}^\gamma (a-x_k)
&\geq  \sqrt{a-1 \over a+1}^{b+2\delta}\sqrt{a^2-1}^{\gamma}\cr 
} 
$$ 
where $\delta = 0$ if $b \in {\bf Z}$, $\delta = 1$ otherwise.
}

\bigskip

 \noindent {\bf Lemma 4.} {\it In the situation of lemma 3, we have
$$\sup_{(x_1, \dots, x_\gamma) \in P} \prod_{k=1}^\gamma (a-x_k) =
(a-{b \over \gamma})^\gamma.
$$
}

\bigskip

We now return to the proof of theorem 2. Since $Pr$ is isogenous to a direct factor of $J_X$ in $J_Y$ by proposition 1, we have
$$\sharp Pr({\bf F}_q) = {\sharp J_Y({\bf F}_q) \over \sharp J_X({\bf F}_q)}. \leqno (4)$$
We need a more precise form of Weil theorem stated in the introduction.
Let $F_C$ be the Frobenius endomorphism acting on the tate module $T_\ell(J_C)$ of the Jacobian $J_C$ of a smooth projective curve $C$ over ${\bf F}_q$. It is well known that $\dim T_\ell(J_C) = 2 g_C$ where $g_C$ denotes the genus of $C$. If $Spec ~F_C$ is the spectrum with multiplicities of this endomorphism, then Weil's theorem asserts that

$$\sharp C({\bf F}_q) = q+1-\sum_{\omega \in Spec ~ F_C} \omega, \leqno (5)$$
$$\sharp J_C(k)= \prod_{\omega \in Spec ~F_C} (1-\omega), \leqno (6)$$
and $\vert \omega \vert = \sqrt q$ for all $\omega \in Spec ~F_C$.

Now, $J_X$ is a $Gal(\overline k/k)$-invariant subvariety of $J_Y$, so that $T_\ell(J_X)$ is an $F_Y$-invariant submodule of $T_\ell(J_Y)$, and $F_X$ is the restriction of $F_Y$ to $T_\ell(J_X)$. Moreover, $\dim T_\ell(J_Y)-\dim T_\ell(J_X)=2g-1-g=g-1$ by Riemann-Hurwitz formula. Hence, there exists some numbers $\theta_1, \dots, \theta_{g-1} \in {\bf R}/2\pi {\bf Z}$, such that
$$Spec ~F_Y= Spec~F_X \cup \left\{\sqrt q \exp (\pm \imath \theta_1),
\cdots, \sqrt q \exp (\pm \imath \theta_{q-1})\right\}.$$
This implies, together with $(5)$ applied to both $X$ and $Y$ :
$$\sharp Y({\bf F}_q)-\sharp X({\bf F}_q) = -2\sqrt q \sum_{k=1}^{g-1} \cos \theta_k, \leqno (7)$$
and together with $(4)$ and $(6)$ applied to both $J_X$ and $J_Y$ :
$$\sharp Pr_\pi(k)=\prod_{k=1}^{g-1}(q+1-2\sqrt q \cos \theta_k). \leqno (8)$$
Notice that we haven thanks to $(7)$
$$-(g-1) \leq \sharp Y({\bf F}_q)-\sharp X({\bf F}_q) \leq g-1.$$

One then deduce part $(\imath)$ of theorem 2 (resp part $(\imath \imath)$) from $(7), (8)$ and lemma 3 (resp. lemma 4) with $C=X$ and $Y$, $a={q+1 \over 2 \sqrt q}$,
$b = -{\sharp Y({\bf F}_q)-\sharp X({\bf F}_q) \over 2 \sqrt q}$, and $\gamma = g-1$. Part $(\imath \imath \imath)$ follows then from the inequalities
$$-d(q+1) \leq -\sharp X({\bf F}_q) \leq \sharp Y({\bf F}_q)-\sharp X({\bf F}_q) \leq \sharp X({\bf F}_q)\leq d(q+1).$$
The proof of the theorem will be complete if we prove both lemma 3 and 4.

\bigskip

 \noindent {\it Proof of lemma 3.} Up to logarithm, we have the calculate the minimum of the function
$$F(x_1, \dots, x_\gamma)= \sum \log(a-x_k).$$
This is a strictly convex function, hence its minimum on the compact convex domain $P$ is reached on an extremal point of $P$, that is a point $x \in P$, such that $[a, b] \subset P$ implies $x \notin ]a, b[$.

Suppose that $x = (x_1, x_2, x_3, \dots, x_\gamma) \in P$ have at least two coordinates, $x_1$ and $x_2$ for simplicity, lying in $]-1, 1[$. Then for $\varepsilon$ small enough, one also have
$(x_1+t, x_2-t, x_3, \dots, x_\gamma) \in P$ for any $t \in ]-\varepsilon, \varepsilon[$, which implies that $x$ is not extremal on $P$. It follows that an extremal point $a = (a_1, \dots, a_\gamma) \in P$ satisfies
$$\cases{ 
a_k=1 &for
$n$ values of $k \in \{1, \dots, \gamma\}$,\cr 
a_k=-1 &for
$m$ values of $k \in \{1, \dots, \gamma\}$, \cr
{\rm eventually~} a_k= \pm \{b\} {\rm ~or } \pm(\{b\}-1) &for $\delta \in \{0, 1\}$ value of $k \in \{1, \dots, \gamma\}$.\cr 
} 
$$
Here, $\{b\}$ denotes the fractionnal part of the real number $b$, so that $\{b\} \in [0, 1[$ and $b-\{b\} \in {\bf Z}$. Hence, we have $\delta = 0$ if $b \in {\bf Z}$, and $\delta = 1$ otherwise.

In case $b \notin {\bf Z}$, that is $\delta = 1$, let us denote by $\beta$ the unique coordinate of the extremal point $a = (a_1, \dots, a_\gamma)$ of $P$, lying in $]-1, 1[$. Up to permutation of the entries, these extremal points have then the shape
$$\cases{ 
(1, \dots, 1, -1, \dots, -1) &if $b \in {\bf Z}$,\cr 
(1, \dots, 1, -1, \dots, -1, \beta) &if $b \notin {\bf Z}$.\cr 
} 
$$
Now, the equations
$$\left\{\matrix{
n+m+\delta &= &\gamma \cr
n-m+\delta \beta &= &b \cr
}\right. \leqno (9)$$
gives easily the values of $n$ and $m$ in terms of the parameters $\gamma, b$ and $\delta$.
We obtain :
$$\eqalign 
{\min \exp F(x_1, \dots, x_\gamma)
&= \exp F(a_1, \dots, a_\gamma) \cr
&= (a- \beta)^\delta(a-1)^n(a+1)^m \cr 
&= (a- \beta)^\delta(a^2-1)^{m+n\over 2}\left({a-1\over a+1}\right)^{n-m\over 2} \cr
&= (a- \beta)^\delta(a^2-1)^{\gamma-\delta\over 2}\left({a-1\over a+1}\right)^{b-\delta \beta \over 2}
} 
$$ 
thanks to $(9)$. But $a-\beta \geq a-1$ and $b-\delta \beta \leq b+\delta$, hence lemma 3.

\bigskip

 \noindent {\it Proof of lemma 4.} We now have to calculate the maximum of the strictly convex function $F$ given in the beginning of the proof of lemma 3. It as a maximum on the compact set $P$.
 Since $F$ is strictly convex on its range of definition, this is also a maximum on the larger set characterized by $\sum x_k=b$.
 By differential calculus, this maximum of $F$ on the whole hyperplane whose equation is $\sum x_k=b$ is reached at the point $x=(x_1, \dots, x_\gamma)$ provided that $grad ~F(x) = (-{1 \over a-x_1}, \dots, -{1\over a-x_\gamma})$ is colinear to
$grad ~(\sum x_k -b) = (1, \dots, 1)$, that is for $x_1= \dots = x_\gamma = {b \over \gamma}$. Since
$\vert b \vert \leq \gamma$ by assumption, this point lies on $P$, hence lemma 4.

\bigskip

 \noindent {\bf III. Remarks} 
 
 \noindent {\bf Remark 1.} Let $\pi : Y \rightarrow X$ be a Galois covering of curves of any degree, and with possible ramifications.
 The Galois group $G$ acts on the Jacobian $J_Y$ and on its Tate module.
 For each irreducible character $\chi$ of  $G$, there is an isotypic sub-abelian variety $P_\chi$ of $J_Y$. 
 For the trivial character $1$, we have $P_1=J_X$.
 One could hope to apply the method of this paper to obtain bounds on the number of rationnal points of these $P_\chi$. Unfortunately, they depend on the dimension of $P_\chi$, which cannot in general be expressed in terms of simple invariants. The point in the present paper is that in the degree two case, there is an unique non-trivial $\chi$, and if moreover the covering is unramified, then one can compute by Riemann-roch the dimension of $J_Y$ in term of the genus $g$ of $X$, so that the dimension of the Prym variety $Pr_\pi = P_\chi$ for the unique non-trivial $\chi$ is known.
 
 However, these dimensions can be computed in case $\pi : Y \rightarrow X$ is a Galois covering of order prime to the characteristic $p$, whose group $G$ have only rational representations (note that this is the case if $G$ is a Wey-group). The reader is refered to [2]. The proofs therein also work in finite characteristic if one works with the Tate modules instead of the cohomology groups $H^0(X, \omega_X)$.

 \bigskip

 \noindent {\bf Remark 2.} Since $Pr$ have dimension $g-1$, we already saw in the introduction that  
$$(q+1-2\sqrt q)^{g-1} \leq \sharp Pr(k) \leq (q+1+2\sqrt q)^{g-1}$$
thanks to Weil theorem. On the other hand, $(6)$ implies
$$-2(g-1) \sqrt q  \leq \sharp Y({\bf F}_q)-\sharp X({\bf F}_q) \leq 2(g-1) \sqrt q.$$
Hence, our bounds in theorem 2 are always ``better" than Weil's one (in the sense that for instance our upper bound is smaller than Weil's one).

Moreover, one can observe that the smaller $\vert \sharp Y({\bf F}_q)-\sharp X({\bf F}_q) \vert $ is  in front of its maximal possible value
$2(g-1) \sqrt q$, the ``sharper" (in the sense that the ratio of our upper bound by our lower bound is smaller) are our bounds in theorem 2. But for fixed $\pi : Y \rightarrow X$ over ${\bf F}_q$ and large $n$, it is hopefull by Tchebotarev that about half of the rational points in $X({\bf F}_{q^n})$ split in $Y$, and half remain inert, so that $\sharp Y({\bf F}_{q^n})-\sharp X({\bf F}_{q^n})$ should be small in front of $2(g-1) \sqrt {q^n}$.
Consequently, it can be expected that  our bound is rather good for large ${q^n}$.

\bigskip

 \noindent {\bf Remark 3.} Of course, one can also give upper and lower bouds for the number of rational points of $J_X$ in the same way.  We obtain :
 
 \medskip
 
\noindent  {\bf Theorem 5.} {\it Let $J_X$ be the jacobian variety of the projective smooth irreducible curve $X$ of genus $g$ defined over ${\bf F}_q$. Then

$$ \left({\sqrt q +1 \over \sqrt q-1}\right)^{{\sharp X({\bf F}_q)-(q+1) \over 2 \sqrt q}-2 \delta} (q-1)^{g}\leq \sharp J_X({\bf F}_q)\leqno (\imath)$$
with 
$\delta = 1$ if ${\sharp X({\bf F}_q)-q-1\over 2 \sqrt q} \notin {\bf Z}$, and $\delta = 0$ otherwise.

$$\sharp J_X({\bf F}_q) \leq \left(q+1 + {\sharp X({\bf F}_q)-q-1\over g} \right)^{g}.\leqno (\imath \imath)$$
 }
 
 \bigskip
 
 \noindent {\it Proof of theorem 5.} It follows from $(5)$, $(6)$, and lemma 3 and 4 with $\gamma = g$, $a = {q+1 \over 2\sqrt q}$ and $b = -{\sharp X({\bf F}_q)-q-1 \over 2 \sqrt q}$.
 
 \bigskip
 
 Let us compare this theorem with Lachaud and Martin-Deschamps one stated in the introduction. Rougthly speaking, their upper and lower bounds, say for fixed $q$ and large $g$,
 grows both like $q^g$.
 The bounds of theorem 5 can be better if $\sharp X({\bf F}_q)$ is sufficiently small.
 However, suppose that $X$ has a degree $d$ map onto the projective line.
 Then $\sharp X({\bf F}_q) \leq d(q+1)$, and the upper bound of proposition 5 implies
 $$\sharp J_X({\bf F}_q) \leq (q+1)^g \left(1 + {d-1\over g} \right)^{g}Ê\leq  \exp(d-1) (q+1)^g,$$
 growing like $(q+1)^g$, which is not as good as $(3~bis)$.

\bigskip

\noindent {\bf Remerciements.} Ce travail a \'et\'e \'elabor\'e dans une grande mesure lors de discutions avec Emmanuel Hallouin et Thierry Henocq. L'auteur tient \`a les remercier vivement.

\bigskip

\centerline{\bf References.}

\medskip

 \noindent [1] A. Beauville, 
Prym varieties and the Schottky problem, {\it Inv. Math.}, 41 (1977), 149-196.

\medskip

 \noindent [2] A. Ksir, 
Dimensions of Prym varieties, {\it Int. Journ. of Math. and Math. Sc.}, 26 (2001), 107-116.

\medskip

 \noindent [3] G. Lachaud, M. Martin-Deschamps, 
Nombre de points des jacobiennes sur un corps fini, {\it Acta Arith.} LVI (1990), 329-340.

\medskip

 \noindent [4] D. Mumford, 
Prym varieties I, in {\it Contribution in Analysis.} LVI (1974), 325-350 ; = Selected paper, (2004) 545-570 Springer.

\medskip

\noindent [5] W. Waterhouse, Abelian Varieties over finite fields, {\it Ann. Sci. \'Ec. Norm. Sup} (4) 2 (1969), 521-560.

\medskip

 \noindent [6] A. Weil, {\it Courbes algebriques et vari\'et\'es ab\'eliennes}, (Hermann, Paris, 1948).

\bigskip

\noindent Marc Perret

\noindent GRIMM

\noindent Universit\'e de Toulouse 2 - Le Mirail

\noindent 5, all\'ees A. Machado

\noindent 31 058 Toulouse - France
 
\noindent perret@univ-tlse2.fr

\end

let $E$ be an elliptic curve over ${\bf F}_q$. There is a degree two map from $E$ to the affine line ${\bf P}^1$, ramified over $4$ geometric points. Let $C$ be a genus $0$ double cover of ${\bf P}^1$ : it is ramified over two geometric points. Assume that these $4+2=6$ points are distinct.
Let $Y=E\times_{{\bf P}^1}C$ be the fiber product of $E$ by $C$ over ${\bf P}^1$ with respect to these degree two maps on ${\bf P}^1$. Then $Y$ is a $\left({\bf Z}/2{\bf Z}\right)^2$-Galois covering of ${\bf P}^1$, with $3$ quadratic intermediate coverings $E, C$ and a certain $X$. Riemann-Hurwitz formula (applied either to the double cover $Y \rightarrow E$ or $Y\rightarrow C$) implies that $Y$ have genus $3$. Moreover, $Y \rightarrow X$ is an unramified double cover. Hence $X$ have genus $g=2$ again by Riemann-Hurwitz applied to $Y \rightarrow X$.
Now, it is easily seen that the pulling back from $E$ to $J_Y$ falls in the set of fixed points of the non-trivial involution of $Y \rightarrow X$. This implies, thanks to dimensions and a connectedness argument, that $E = Pr_{Y \rightarrow X}$ is a prymian variety.